\documentclass[11pt]{elsarticle}
\usepackage{amssymb, color}
\usepackage{amsmath}
 \usepackage{pifont}
 \usepackage[utf8]{inputenc}

 \usepackage[utf8]{inputenc}
\usepackage[T1]{fontenc}
\usepackage[french]{babel}
\usepackage{enumitem}

\makeatletter
\def\LaTeX{\leavevmode L\raise.42ex
    \hbox{\kern-.3em\size{\sf@size}{0pt}\selectfont A}\kern-.15em\TeX}
\makeatother

\newcommand{\BibTeX}{{\rm B\kern-.05em{\sc
          i\kern-.025emb}\kern-.08em\TeX}}

\makeatletter
\def\@currentlabel{2.1}\label{e:dispaa}
\def\@currentlabel{2.21}\label{e:dispau}
\def\@currentlabel{2.22}\label{e:dispav}
\def\@currentlabel{2.23}\label{e:dispaw}
\def\@currentlabel{2.24}\label{e:dispax}
\def\theequation{\thesection.\@arabic\c@equation}
\makeatother

\renewcommand{\theequation}{\arabic{section}.\arabic{equation}}
\def \D{\Delta}

\def \l{\lambda}

\def \p{\partial}

\newcommand{\R}{\mathbb R}

\newcommand{\N}{\mathbb N}

\def \D{\Delta}

\newtheorem{thm}{Theorem} [section]
\newtheorem{lem}{Lemma} [section]
\newtheorem{prop}{Proposition} [section]
\newtheorem{cor}{Corollary} [section]

\newtheorem{rem}{Remark}[section]

\renewcommand{\theequation}{\thesection.\arabic{equation}}
\renewcommand{\thesection}{\arabic{section}}
\renewcommand{\theequation}{\thesection.\arabic{equation}}
\let\ssection=\section\renewcommand{\section}{\setcounter{equation}{0}\ssection}
\begin{document}
\begin{frontmatter}

\title{Liouville theorems for stable solutions of the weighted Lane-Emden system}

\author[hh]{Hatem hajlaoui}
\ead{hajlouihatem@gmail.com}
\author[ah]{Abdellaziz Harrabi}
\ead{abdellaziz.harrabi@yahoo.fr}
\author[fd]{Foued Mtiri\corref{cor1}}
\ead{mtirifoued@yahoo.fr}
\address[hh]{Institut Pr\'{e}paratoire aux Etudes d'Ing\'{e}nieurs. Universit\'e de Kairouan, Tunisia.}
\address[ah]{Institut Sup\'erieur des Math\'ematiques Appliqu\'ees et de l'Informatique. Universit\'e de Kairouan, Tunisia.}
\address[fd]{Institut \'{E}lie Cartan de Lorraine, IECL,  UMR 7502, Universit\'{e} de Lorraine, France.}
\begin{abstract} We examine the general weighted Lane-Emden system
\begin{align*}
 -\Delta u = \rho(x)v^p,\quad -\Delta v= \rho(x)u^\theta, \quad u,v>0\quad \mbox{in }\;\mathbb{R}^N
\end{align*}
 where $1<p\leq\theta$ and
 $\rho: \mathbb{R}^N\rightarrow \mathbb{R}$ is a radial continuous function satisfying
 $\rho(x)\geq A(1+|x|^2)^{\frac{\alpha}{2}}$ in $\mathbb{R}^N$ for some $\alpha\geq 0$ and $A>0$. We
 prove some Liouville type results for stable solution and improve the previous works \cite{co, Fa, HU}. In particular, we establish a new comparison property
 (see Proposition \ref{p1} below) which is crucial to handle the case $1 < p \leq \frac{4}{3}$. Our results can be applied also to the weighted Lane-Emden equation
$-\Delta u = \rho(x)u^p$ in $\mathbb{R}^N$.
 \end{abstract}
\end{frontmatter}
\section{Introduction}
\setcounter{equation}{0}
We consider the following weighted Lane-Emden system
\begin{align}\label{1.1}
-\Delta u = \rho(x)v^p, \quad-\Delta v= \rho(x)u^\theta,\quad u,v>0\quad\mbox{in }\; \mathbb{R}^N
\end{align}
where $1<p\leq\theta$ and $\rho: \mathbb{R}^N\rightarrow \mathbb{R}$ is a radial continuous function satisfying the following assumption:
\begin{enumerate}
 \item [$(\star)$] There exists $\alpha \geq 0$ and $A>0$ such that $\rho(x)\geq A\rho_0(x)$ in $\R^N$ where $\rho_0
 :=(1+|x|^2)^{\frac{\alpha}{2}}$.
\end{enumerate}
Remark that under the scaling transformation $u =\gamma^{\frac{1}{\theta+1}}\widetilde{u}$, $v = \widetilde{v}$ with $\gamma > 0$,
the following system
\begin{align*}
-\Delta \widetilde{u} = \widetilde{\rho}(x)\widetilde{v}^p, \quad-\Delta \widetilde{v}= \gamma\widetilde{\rho}(x)\widetilde{u}^\theta,
\quad \widetilde{u},\widetilde{v}>0 \quad \mbox{in}\;\; \mathbb{R}^N
\end{align*}
is equivalent to \eqref{1.1} with $\rho = \gamma^{\frac{1}{\theta+1}}\widetilde{\rho}.$

\medskip
To define the notion of stability, we consider a general system given by
\begin{align}\label{1.2}
 -\Delta u = f(x,v),\quad -\Delta v= g(x,u),\quad \mbox{in }\; \mathbb{R}^N
\end{align}
with $f,g\in C^1(\mathbb{R}^{N+1},\mathbb{R})$ satisfying $f_s:=\frac{\partial f(x,s)}{\partial s},
g_s=\frac{\partial f(x,s)}{\partial s}\geq 0$ in $\R.$ A smooth solution $(u,v)$ of \eqref{1.2} is said stable if there
exist positive smooth functions $\xi$, $\zeta$ verifying
\begin{align*}
 -\Delta \xi = f_{v}(x,v)\zeta, \quad -\Delta \zeta = g_{u}(x,u)\xi,\quad \mbox{in }\, \mathbb{R}^N.
\end{align*}
This definition is motivated by \cite{Mo, Fa, co}. In this paper, we prove the following classification result:
\begin{thm}
\label{main1} Suppose that $\rho$ satisfies $(\star)$ and let $x_0$ be the largest root of the polynomial $H(x) = $
\begin{align}\label{H}
x^4 -
\frac{16p\theta(p+1)(\theta+1)}{(p\theta-1)^2}x^2 +
\frac{16p\theta(p+1)(\theta+1)(p+\theta+2)}{(p\theta-1)^3}x
-\frac{16p\theta(p+1)^2(\theta+1)^2}{(p\theta-1)^4}.
\end{align}
 \begin{enumerate}[label=\roman*)]
\item
 If $\frac{4}{3}< p \leq \theta$ then \eqref{1.1} has no stable classical solution if $N<2+(2+\alpha)x_0.$
 In particular, if $N \leq 10+4\alpha,$ then \eqref{1.1} has no classical stable
solution for all $\frac{4}{3}<  p \leq \theta$.
\item If $1<p\leq \min(\frac{4}{3}, \theta)$, then \eqref{1.1} has no bounded classical stable solution, if
$$N < 2 + \left[\frac{p}{2}+\frac{(2-p)(p \theta -1)}{(\theta+p-2)(\theta+1)}\right](\alpha+2)x_0.$$
\end{enumerate}
Therefore, if $N \leq 6+2\alpha,$ the system \eqref{1.1} has no bounded classical stable
solution for all $\theta\geq p > 1.$
\end{thm}
As a consequence of Theorem \ref{main1}, we obtain the following classification result for stable solution of the Lane-Emden equation
\begin{align}\label{D}
 -\Delta u = \rho(x)u^p, \quad u>0 \quad \mbox{in }\;\mathbb{R}^N.
\end{align}
\begin{cor}\label{main2} Suppose that $\rho$ satisfies $(\star)$ and let $p>1.$
\begin{enumerate}[label=\roman*)]
 \item
 If $\frac{4}{3}< p$ then \eqref{D} has no stable classical solution if
\begin{align}\label{N}
N<2+\frac{2(2+\alpha)}{p-1}\left(p+\sqrt{p^2-p}\right).
 \end{align}
 In particular, if $N \leq 10+4\alpha,$ then \eqref{D} has no stable classical
solution for all $\frac{4}{3}< p$.
\item If $1<p\leq \frac{4}{3} $,
\eqref{D} has no bounded stable classical solution for $N$ verifying \eqref{N}.
\end{enumerate}
Therefore, there is no bounded stable classical solution of \eqref{D} for all $p > 1$ if $N \leq 10+4\alpha.$
\end{cor}

Recalling that for the autonomous case, i.e. when $\rho\equiv 1,$ the stable solutions of the corresponding Lane-Emden equation and
system, or the biharmonic equation (corresponding to $p = 1$) have been widely studied by many authors. See for instance \cite{Far, wy, co, hhy, cdg, ddww} and the references there in.

\medskip
For the second order Lane-Emden equation ($p > 1$)
\begin{align}\label{1.4}
 -\Delta u = |u|^{p-1}u\quad \mbox{in }\;\mathbb{R}^N,
\end{align}
Farina classified completely in \cite{Far} all finite Morse index classical solutions for $1<p<p_{JL},$ where $p_{JL}$ stands for the Joseph-Lundgren exponent \cite{jl} (see also \cite{gnw}).
More precisely, the equation \eqref{1.4} admits nontrivial classical solutions with finite Morse index if
and only if $N\geq3$, $p=\frac{N+2}{N-2}$ or $N\geq11$ and $p\geq p_{JL}.$ For the biharmornic equation ($p > 1$)
\begin{align}\label{1.5}
 \Delta^2 u = |u|^{p-1}u\quad \mbox{in }\;\mathbb{R}^N,
\end{align}
 D\'avila-Dupaigne-Wang-Wei \cite{ddww} recently gave a complete classification of finite Morse index solutions.
They derived a monotonicity formula for the solutions of \eqref{1.5} and reduced the problem to the nonexistence of stable homogeneous
solutions.

\medskip
It is worthy to mention that Chen-Dupaigne-Chergu \cite{cdg} proved an optimal Liouville type result for the radial stable solutions of \eqref{1.1}
for $\theta \geq p > 1$ and $\rho \equiv 1$.

\medskip
For the weighted equation or system with positive weights, the Liouville type results are less understood.
 \begin{itemize}
   \item Using Farina's approach, Fazly proved the nonexistence of classical stable solutions of \eqref{D} for $\rho=\rho_0,$ $N$ satisfying \eqref{N} and $p\geq 2.$ See Theorem 2.3 in \cite{Fa}.
   \item Using also Farina's approach, Cowan-Fazly \cite{cf} established a Liouville type result for classical stable sub-solutions
   of \eqref{D} for $N$ satisfying \eqref{N}, $p>1$ and
        \begin{align}\label{cc}
         \displaystyle{\lim_{|x|\rightarrow \infty}}\frac{\rho(x)}{\rho_0(x)}=C\in (0,\infty).
        \end{align}
        See Theorem 1.3-(3) with $\alpha = 0$ in \cite{cf}.
   \item Adopting the new approach in \cite{cg}, Hu proved the following Liouville theorem for classical stable solutions of \eqref{1.1} for $\rho=\rho_0$ and $\theta\geq p\geq 2$ or $\theta=p>\frac{4}{3},$ obtaining a direct extension of Theorem 1 in \cite{co} for $\rho\equiv 1$.
       More precisely, let $t_0^+$ and $t_0^-$ be the quantities used in \cite{co}:
\begin{align*}
t_{0}^{\pm} = \sqrt{\frac{p\theta(p+1)}{\theta+1}}\pm\sqrt{\frac{p\theta(p+1)}{\theta+1}-\sqrt{\frac{p\theta(p+1)}{\theta+1}}},
\end{align*}
\end{itemize}
Hu proved in \cite{HU}:
\begingroup
\setcounter{thm}{0} %assign desired value to theorem counter
\renewcommand\thethm{\Alph{thm}}% locally redefine the representation of the theorem counter
\begin{thm}
\label{HU} Suppose that $\rho=\rho_0$ with $\alpha\geq0.$
 \begin{enumerate}[label=\roman*)]
\item
 If $2t_0^-<p\leq \theta$ and $N$ satisfies
$$N<2+ \frac{2(2+\alpha)(\theta+1)}{p\theta-1}t_0^+,$$
then there is no classical stable solution of \eqref{1.1}. In particular there is no classical stable solution of \eqref{1.1} for any
$2\leq p\leq \theta$ and $N\leq 10+4\alpha.$
\item If $p>\frac{4}{3} $ and $N$ satisfies \eqref{N}, then there is no classical stable solution of \eqref{D}.
\end{enumerate}
\end{thm}
\endgroup

\begin{rem}\label{rnew} It is known that for $1<p\leq \theta,$ there hold
$t_{0}^-<1<t_{0}^+,$ $t_{0}^-$ is decreasing and $t_{0}^+$ is increasing in $z:=\frac{p\theta(p+1)}{\theta+1}.$ Moreover,
$\lim_{z\rightarrow \infty}t_0^-= \frac{1}{2}$ and $\lim_{z\rightarrow \infty}t_0^+= 1.$
\end{rem}

\begin{rem}
\label{rnew1}
We have $2t_0^-<p$ if $p>\frac{4}{3}$. Indeed, if $p>\frac{4}{3}$ then $\theta\geq p>\frac{4}{3}$ and $z>\frac{16}{9}.$
Since $f(z):=\sqrt{z}-\sqrt{z-\sqrt{z}} $ is decreasing in $z,$ there holds $2t_0^-=2f(z)<2f(\frac{16}{9})=\frac{4}{3} < p.$
\end{rem}

Using the above remark, we see that Theorem \ref{HU} (hence Theorem 1 in \cite{co}) can be extended immediately for
$\frac{4}{3}<p\leq \theta$.
\begin{itemize}
  \item We can show that $2t_0^+\frac{\theta+1}{p\theta-1}<x_0$ for any $1<p\leq \theta$ (see Lemma \ref{l.2.4} below), where
  $x_0$ is the largest root of the polynomial $H$ given by \eqref{H}. So Theorem \ref{main1} improves the bound given in Theorem \ref{HU}.
  \item In Theorem \ref{main1} and Corollary \ref{main2}, we prove classification results for \eqref{1.1} and \eqref{D} with
  $\rho$ satisfying $(\star)$, so without the restriction $\rho=\rho_0$ in Theorem \ref{HU}; or the condition \eqref{cc} used in \cite{cf}.
  \item Our approach permits to prove a Liouville type result for $\theta \geq p> 1.$ To the best of our knowledge, no general Liouville type
  result was known for stable solution of \eqref{1.1} with positive weight for $1< p \leq \frac{4}{3}.$
\end{itemize}

To prove Theorem \ref{main1}, we will use the following Souplet type estimate \cite{S}. Its proof is the same as for Lemma 2.3 in
\cite{HU} where we replace just $\rho_0$ by $\rho$, so we omit the details.
\begin{lem}\label{Soup} Let $\theta \geq p > 1$ and $\rho$ satisfy $(\star).$ Then any classical solution of \eqref{1.1}
verifies
\begin{align}
\label{estS}
u^{\theta + 1} \leq \frac{\theta+1}{p+1}v^{p+1} \;\; \mbox{in}\;\;\mathbb{R}^N.
\end{align}
\end{lem}
However, to handle the case $1 < p \leq \frac{4}{3}$, we need the following new comparison property between $u$ and $v$. It is
somehow an inverse version of Souplet's estimate \eqref{estS}, and has its own interest.
\begin{prop}\label{p1}
Let $\theta\geq p>1$ and suppose that $\rho$ satisfies $(\star).$ Let $(u,v)$ be a classical solution of \eqref{1.1} and assume
that $u$ is bounded, then
\begin{align*}
v\leq \|u\|_\infty^\frac{\theta-p}{p+1}u.
 \end{align*}
\end{prop}

\medskip
Our paper is organized as follows. In section 2, we prove some preliminaries results, in particular we give the proof
of Proposition \ref{p1}. The proofs of Theorem \ref{main1} and Corollary \ref{main2} are given in section 3.

\section{Preliminaries}
\setcounter{equation}{0}
In order to prove our results, we need some technical lemmas. In the following, $C$
denotes always a generic positive constant independent on $(u,v)$, which could be changed from one line to another.
The ball of center $0$ and radius $r > 0$ will be denoted by $B_r$.

\subsection{Comparison property}
In this subsection, we give the proofs of Proposition \ref{p1}. First, we can adapt the proof of Lemma 2.1 in \cite{Fa} (which was inspired by the previous works \cite{sz, MPO}), to obtain the following integral estimates for all classical solutions of \eqref{1.1}.
\begin{lem}\label{l.2.1}
Let $p\geq1,\,\theta>1$ and suppose that $\rho$ satisfies $(\star).$ For any classical solution $(u, v)$ of \eqref{1.1}  there exists $C > 0$ such that for any $R \geq 1$, there holds
\begin{align*}
   \int_{B_{R}}\rho(x)v^{p} dx\leq C R^{N-\frac{2(\theta+1)p}{p\theta-1}-\frac{(p+1)\alpha}{p\theta-1}},\quad
   \int_{B_{R}}\rho(x)u^{\theta} dx \leq C R^{N-\frac{2(p+1)\theta}{p\theta-1}-\frac{(\theta+1)\alpha}{p\theta-1}}.
\end{align*}
\end{lem}
\textbf{Proof.}
  Let $\varphi_0 \in C_c^\infty(B_2)$ be a cut-off function verifying $0 \leq \varphi_0 \leq 1$, $\varphi_0=1$ for $|x|<1$. Take
  $\psi := \varphi_0(R^{-1}x)$ for $R \geq 1.$ Multiplying the equation $-\D u=\rho(x)v^{p}$ by $\psi ^{m}$ and integrating by parts,
  there holds then
\begin{align*}
   \int_{\R^N}\rho(x)v^{p}\psi^{m} dx = -\int_{\R^N}u\Delta(\psi^{m}) dx \leq \frac{C}{R^{2}}\int_{B_{2R}\setminus B_{R}} u \psi^{m-2} dx.
\end{align*}
By H\"older's inequality, we get
\begin{align*}
  \int_{\R^N}\rho(x)v^{p}\psi^{m}dx  \leq& \; \frac{C}{R^{2}}\left(\int_{B_{2R}\setminus B_{R}} \rho(x)^{-\frac{\theta'}{\theta}}dx\right)^{\frac{1}{\theta'}}\left(\int_{B_{2R}\setminus B_{R}}\rho(x) u^{\theta} \psi^{(m-2)\theta}dx\right)^{\frac{1}{\theta}},
\end{align*}
where $\frac{1}{\theta}+\frac{1}{\theta'}=1.$ From $(\star)$ we deduce that for $R \geq 1$,
\begin{align*}
  \int_{\R^N}\rho(x)v^{p}\psi^{m}dx  \leq& \; C R^{\frac{N}{\theta'} -\frac{\alpha}{\theta}-2}\left(\int_{B_{2R}\setminus B_{R}}\rho(x) u^{\theta} \psi^{(m-2)\theta}dx\right)^{\frac{1}{\theta}}.
\end{align*}
 Similarly, using $-\Delta v = \rho(x)u^p$, we obtain, for $k \geq 2$,
\begin{align*}
  \int_{\R^N}\rho(x)u^{\theta}\psi^{k} dx \leq& \; C R^{\frac{N}{p'}-\frac{\alpha}{p} -2}\left(\int_{B_{2R}\setminus B_{R}}\rho(x) v^{p} \psi^{(k-2)p}dx\right)^{\frac{1}{p}},
\end{align*}
where $\frac{1}{p}+\frac{1}{p'}=1.$ Take now $k$ and $m$ large verifying $ m \leq (k-2)p$ and $k \leq (m-2)\theta.$ Combining the two
above inequalities, we get
\begin{align*}
  \int_{\R^N}\rho(x)v^{p}\psi^{m}dx & \leq C R^{\frac{N}{\theta'} -\frac{\alpha}{\theta}-2} R^{\big(\frac{N}{p'} -\frac{\alpha}{p}-2\big)\frac{1}{\theta}}\left(\int_{B_{2R}\setminus B_{R}}\rho(x) v^{p} \psi^{(k-2)p}dx\right)^{\frac{1}{p\theta}}\\
  & \leq C R^{N - \frac{N}{p\theta} -\frac{\alpha(p+1)}{p\theta}- \frac{2(\theta+1)}{\theta}}\left(\int_{\R^N}\rho(x) v^{p} \psi^m dx\right)^{\frac{1}{p\theta}}.
\end{align*}
Hence
\begin{align*}
  \int_{B_R}\rho(x)v^{p} dx\leq \int_{\R^N}\rho(x)v^{p}\psi^{m}dx \leq C R^{N-\frac{2(\theta+1)p}{p\theta-1}-\frac{(p+1)\alpha}{p\theta-1}}.
\end{align*}
Similarly, we obtain the estimate for $u$.\qed

\medskip
Now we are in position to prove the inverse comparison property.

\noindent
{\bf Proof of Proposition \ref{p1}.} Let $w= v- \l u,\;$ where $\l = \|u\|_{\infty}^\frac{\theta-p}{p+1}$. We have, as $\theta \geq p$,
\begin{align*}
\D w = \rho(x)\left(\l v^p- u^\theta\right) =\rho(x)\left[ \l v^p- \left(\frac{u}{ \|u\|_{\infty}}\right)^{\theta} \|u\|_{\infty}^{\theta}\right] &\geq \rho(x)\left[\l v^p- \left(\frac{u}{ \|u\|_{\infty}}\right)^p \|u\|_{\infty}^{\theta}\right]\\
&=\rho(x) \|u\|_{\infty}^{\theta - p} \left(\frac{\l v^p}{\|u\|_{\infty}^{\theta - p}} - u^p\right)\\
&= \rho(x)\|u\|_{\infty}^{\theta - p} \left(\l^{-p}v^p - u^p\right).
\end{align*}
It follows that $\D w \geq0$ in the set $\{w\geq 0\}$. Consider $w_+ := \max(w, 0).$ Next, we split the proof into two cases.

\medskip
\textit{Case 1: $p\geq 2$.} For any $R > 0$, there holds
\begin{align}
\label{2.2}
\begin{split}
(p-1) \int_{B_R} w_+^{p-2}|\nabla w_+|^2 dx & = -\int_{B_R} w_+^{p-1}\Delta w dx + \int_{\p B_R}w_+^{p-1} \frac{\p w}{\p\nu}
d\sigma\\ & \leq \int_{\p B_R}w_+^{p-1} \frac{\p w}{\p\nu} d\sigma\\ & = \frac{R^{N-1}}{2}f'(R)
\end{split}
\end{align}
where
$$f(R) := \int_{S^{N-1}} w_+^p(R\sigma)d\sigma \leq \int_{S^{N-1}} v^p(R\sigma)d\sigma =:g(R).$$
Hereafter, $S^{N-1}$ denotes by the unit sphere in $\mathbb{R}^N.$ By Lemma \ref{l.2.1}, we derive that
\begin{align*}
\int_{0}^R r^{N-1}\int_{S^{N-1}}\rho(r\sigma) v^p(r\sigma)d\sigma dr & = \int_{0}^R r^{N-1}\rho(r)\int_{S^{N-1}} v^p(r\sigma)d\sigma dr \\
 \leq&\; C R^{N-\frac{2(\theta+1)p}{p\theta-1}-\frac{(p+1)\alpha}{p\theta-1}} = o(R^N) \quad \mbox{as } R \to \infty.
\end{align*}
Using $(\star),$ there holds
\begin{align*}
\int_{0}^R r^{N-1+\alpha}g(r)  dr = o(R^N) \quad \mbox{as } R \to \infty.
\end{align*}
This implies that $\liminf_{r\to \infty} g(r) = 0,$ hence $\liminf_{r\to \infty} f(r) = 0.$ Consequently, there exist
$R_i \to \infty$ such that $f'(R_i) \leq 0$. Take \eqref{2.2} with $R = R_i$ and let $i\rightarrow \infty$, we conclude that $w_+$ is constant in $\mathbb{R}^N.$ If $w \equiv C>0$ then $v\geq C>0$ in $\mathbb{R}^N,$ which contradicts Lemma \ref{l.2.1}. Hence $w_+ \equiv 0$ in $\mathbb{R}^N$, i.e.~$v - \l u \leq 0$ in  $\mathbb{R}^N.$

\medskip
\textit{Case 2: $1 <p <2.$} For any $R>0$ and $\epsilon >0,$ we have
\begin{align*}
(p-1)\int_{ B_R} (\epsilon + w _+)^{p-2}|\nabla w_+|^2 dx & = -\int_{ B_R} (\epsilon + w _+)^{p-1}\Delta w dx
 + \int_{\p B_R}(\epsilon + w _+)^{p-1} \frac{\p w}{\p\nu} d\sigma\\
 & \leq \int_{\p B_R}(\epsilon + w _+)^{p-1} \frac{\p w}{\p\nu} d\sigma.
\end{align*}
Letting $\epsilon \rightarrow 0$ (passing to limit in the l.h.s.~via monotone convergence and use the dominated convergence
on the r.h.s.), we get always the estimate \eqref{2.2}, which will lead to the same conclusion: $w_+ \equiv 0$ in $\mathbb{R}^N$.
\qed

\medskip
\subsection{Consequence of stability}

 With the ideas in \cite{cg, dfs}, we can proceed similarly as the proof of Lemma 2.1 in \cite{HU} and claim
 \begin{lem}\label{l.2.2}
If $(u,v)$ is a nonnegative classical stable solution of \eqref{1.1}, then
\begin{equation}
\label{1.3}  \sqrt{p\theta}\int_{\R^N}\rho(x)
u^{\frac{\theta-1}{2}}v^{\frac{p-1}{2}}\phi^2dx  \leq
\int_{\R^N}|\nabla\phi|^2dx , \quad \forall\; \phi \in C_c^1(\R^N).
\end{equation}
\end{lem}

The following Lemma is a consequence of the stability inequality \eqref{1.3} and Proposition \ref{p1}. It plays a
crucial role to handle the case $1<p\leq \frac{4}{3}.$ Here we use also some ideas coming from \cite{wy, hhy}.
\begin{lem}
\label{lemnew}
Let $(u,v)$ be a stable solution to \eqref{1.1} with $1< p \leq \min(\frac{4}{3}, \theta)$. Assume that $u$ is bounded and $\rho$ satisfies $(\star)$, there holds
\begin{align}
\label{2.3}
\int_{B_R} \rho(x)v^2 dx \leq CR^{N-\frac{2(\theta+1)p}{p\theta-1}-\frac{(p+1)\alpha}{p\theta-1} - \frac{2(2+\alpha)(2 - p)}{\theta+p-2}}, \quad \forall\; R > 0.
\end{align}
\end{lem}

\noindent{\bf Proof.} Let $(u,v)$ be a stable solution of \eqref{1.1}, where $u$ is bounded. Take $\eta \in C_c^\infty(\R^N)$. Multiplying $-\D v= \rho(x)u^{\theta}$ by $v \eta ^2$ and integrating by parts, there holds
\begin{align*}
  \int_{\mathbb{R}^N} |\nabla v|^2 \eta^2dx= \int_{\mathbb{R}^N}\rho(x) u^\theta v \eta^2dx + \frac{1}{2}\int_{\mathbb{R}^N}v^2 \D (\eta^2)dx.
\end{align*}
Using Lemma \ref{Soup}, we get
\begin{align*}
  \int_{\mathbb{R}^N} |\nabla v|^2 \eta^2dx\leq  \sqrt{\frac{\theta+1}{p+1}}\int_{\mathbb{R}^N}\rho(x) u^{\frac{\theta-1}{2}} v^{\frac{p+1}{2}}v \eta^2dx + \frac{1}{2}\int_{\mathbb{R}^N}v^2 \D (\eta^2)dx.
\end{align*}
Set $\phi = v \eta$ in \eqref{1.3} and integrating by parts, we deduce that
\begin{align*}
 \sqrt{p\theta}\int_{\mathbb{R}^N} \rho(x)u^{\frac{\theta-1}{2}} v^{\frac{p-1}{2}}v^2 \eta^2dx \leq  \int_{\mathbb{R}^N} |\nabla v|^2 \eta^2dx +  \int_{\mathbb{R}^N}v^2 |\nabla \eta |^2dx -\frac{1}{2}\int_{\mathbb{R}^N}v^2 \D (\eta^2)dx.
\end{align*}
Combining the two last inequalities, we obtain
 \begin{align*}
 \left(\sqrt{p \theta}-\sqrt{\frac{\theta+1}{p+1}}\right)\int_{\mathbb{R}^N}\rho(x) u^{\frac{\theta-1}{2}} v^{\frac{p+3}{2}} \eta^2dx \leq   \int_{\mathbb{R}^N}v^2 |\nabla \eta |^2dx.
\end{align*}
Using Proposition \ref{p1}, there exists a positive constant $C$ such that
  \begin{align*}
 \int_{\mathbb{R}^N}\rho(x)v^{\frac{\theta+p+2}{2}}\eta^2dx  \leq C \int_{\mathbb{R}^N}v^2 |\nabla \eta |^2dx.
\end{align*}
 Take $\varphi_0$ a cut-off function in $C_c^\infty(B_2)$ verifying $0 \leq \varphi_0 \leq 1$, $\varphi_0=1$ for $|x|<1$. Let $\eta = \varphi^m$ with $\varphi := \varphi_0(R^{-1}x)$ for $R > 0$, we arrive at
 \begin{equation*}
 \int_{\mathbb{R}^N}\rho(x)v^{\frac{\theta+p+2}{2}}\varphi^{2m}dx \leq \frac{C}{R^2}\int_{B_{2R}\setminus B_R} v^2\varphi^{2m-2}dx.
 \end{equation*}
 Using $(\star),$ there holds
 \begin{align}\label{2.13new}
  \int_{\R^N} \rho(x)v^\frac{\theta+p+2}{2}\varphi^{2m}dx
  \leq \frac{C}{R^{2+\alpha}}
 \int_{B_{2R}\backslash B_R} \rho(x)v^2\varphi^{2m-2}dx
  \leq \frac{C}{R^{2+\alpha}}
 \int_{\R^N} \rho(x)v^2\varphi^{2m-2}dx.
 \end{align}
   Denote
 \begin{align*}
J_1 := \int_{\mathbb{R}^N}\rho(x)v^{\frac{\theta+p+2}{2}}\varphi^{2m}dx , \quad  J_2 := \int_{\mathbb{R}^N}\rho(x) v^2\varphi^{2m-2}dx.
 \end{align*}
Remark that $p < 2 < \frac{\theta + p + 2}{2}$ for $1 < p \leq \frac{4}{3}$ and $\theta\geq p$. A direct calculation yields
 \begin{align*}
2 = p\lambda + \frac{\theta + p + 2}{2}(1 - \lambda) \quad \mbox{with } \lambda = \frac{\theta + p - 2}{\theta + 2 - p} \in (0, 1).
 \end{align*}
Take $m$ large such that $m\l > 1.$ By H\"older's inequality, Lemma \ref{l.2.1}  and \eqref{2.13new}, we get
  \begin{align*}
J_2 \leq J_1^{1 - \l} \left(\int_{\R^N}\rho(x)v^p\varphi^{2m\l - 2}dx\right)^\l & \leq \left(\frac{CJ_2}{R^{2+\alpha}}\right)^{1 - \l} \left(\int_{B_{2R}}\rho(x)v^pdx\right)^\l \\ & \leq C'J_2^{1 - \l}R^{-(2+\alpha)(1 - \l)}\left(R^{N-\frac{2(\theta+1)p}{p\theta-1}-\frac{(p+1)\alpha}{p\theta-1}}\right)^\l,
 \end{align*}
which implies
\begin{align*}
J_2 \leq  CR^{N-\frac{2(\theta+1)p}{p\theta-1}-\frac{(p+1)\alpha}{p\theta-1} - \frac{2(2+\alpha)(2 - p)}{\theta+p-2}},
 \end{align*}
so we are done. \qed

\medskip
\subsection{Property of the polynomial $H$}
Consider the polynomial $H$ given by \eqref{H}. Performing the change of variables $x=\frac{\theta+1}{p\theta-1}s$,
a direct computation yields
$$H(x)=\left(\frac{\theta+1}{p\theta-1}\right)^4L(s)$$
where
\begin{align}\label{L}
 L(s):=s^4-\frac{16p\theta(p+1)}{\theta+1}s^2+\frac{16p\theta(p+1)(p+\theta+2)}{(\theta+1)^2}s-\frac{16p\theta(p+1)^2}{(\theta+1)^2}.
\end{align}
 Hence $ H(x)<0$ if and only if $L(s)<0$.
\begin{lem}\label{l.2.4}
Let $1< p \leq \theta$, then $L(2) < 0$ and $L$ has a unique root $s_0$ in $(2, \infty)$ and $2t_0^+ < s_0$. Moreover, if $p > \frac{4}{3}$,
then $L(p) < 0$ and $s_0$ is the unique root of $L$ in $(p, \infty)$.
\end{lem}

\noindent{\bf Proof.} Using $1< p \leq \theta$,
\begin{align*}
L(2) & = 16 -\frac{64p\theta(p+1)}{(\theta+1)} +\frac{32p\theta(p+1)(p+\theta+2)}{(\theta+1)^2} -\frac{16p\theta(p+1)^2}{(\theta+1)^2}\\
& = 16 -\frac{64p\theta(p+1)}{(\theta+1)} + \frac{32p\theta(p+1)}{(\theta+1)} + \frac{32p\theta(p+1)^2}{(\theta+1)^2}
-\frac{16p\theta(p+1)^2}{(\theta+1)^2}\\
& = 16 -\frac{32p\theta(p+1)}{(\theta+1)} + \frac{16p\theta(p+1)^2}{(\theta+1)^2}\\
& \leq 16 - \frac{32p\theta(p+1)}{(\theta+1)} + \frac{16p\theta(p+1)}{(\theta+1)}\\
& = 16 \frac{(1 - p^2)\theta + (1 - p\theta)}{(\theta+1)} < 0.
\end{align*}
Very similarly, we can check that
\begin{align*}
L'(2) \leq 32 - \frac{32p\theta(p+1)}{(\theta+1)} < 0.
\end{align*}
Furthermore, we have
$$L''(s)=12s^2-\frac{32p\theta(p+1)}{\theta+1},$$ then
$L''$ can change at most once the sign from negative to positive for $s \geq 2$. As
$\lim_{s\rightarrow\infty}L(s)= \infty$, it's clear that $L$ admits a unique root in $(2,\infty).$ Moreover, we can check that
$$L(2t_0^+)=\frac{16p\theta(p+1)(\theta-p)}{(\theta+1)^2}(1-2t_0^+) < 0.$$
Hence, there holds $2t_0^+ < s_0.$

\medskip
Now we consider $L(p)$. Rewrite
\begin{align*}
L(s) = s^4 - 16\frac{p\theta(p+1)}{\theta + 1}\left(s^2 - \frac{p+\theta+2}{\theta+1}s + \frac{p+1}{\theta+1}\right).
\end{align*}
For $s > 1$, we see that
\begin{align*}
\left(s^2 - \frac{p+\theta+2}{\theta+1}s + \frac{p+1}{\theta+1}\right)_\theta' = \frac{p+1}{(\theta+1)^2}(s - 1) > 0,
\end{align*}
Then for $s > 1$, as $\theta \geq p > 1$, there holds
\begin{align*}
s^2 - \frac{p+\theta+2}{\theta+1}s + \frac{p+1}{\theta+1} > s^2 - 2s + 1 = (s-1)^2 \quad\mbox{and} \quad \frac{p\theta(p+1)}{\theta + 1} \geq p^2.
\end{align*}
Finally, we get (for $p > 1$)
\begin{align*}
L(p) < p^4 - 16p^2(p-1)^2=p^2(5p-4)(4-3p)
\end{align*}
and
\begin{align*}
L'(p) = 4p^3 - 16\frac{p\theta(p+1)}{\theta + 1}\left(2p - \frac{p+\theta+2}{\theta+1}\right) < 4p^3 - 16p^2(2p - 2) = 4p^2(8 - 7p),
\end{align*}
We check readily that for $p > \frac{4}{3}$, $L(p) < 0$ and $L'(p) < 0$, so we can conclude as above.\qed

\section{Proofs of Theorem  \ref{main1} and Corollary \ref{main2}.}
\setcounter{equation}{0}
The following lemma plays an important role in dealing with Theorems \ref{main1} and Corollary \ref{main2}, where we use some ideas from \cite{hhy}. Here and in the following, we define $R_{k} = 2^{k}R$ for all $R > 0$ and integers $k\geq 1$.
\begin{lem}
\label{newl}Suppose that $\rho$ satisfies $(\star)$ and let $(u,v)$ be a stable solution of \eqref{1.1}. Then for any $s > \frac{p+1}{2}$ verifying
$L(s) < 0$, there exists $C <\infty$ such that
\begin{equation*}
\int_{B_R}\rho(x)u^{\theta}v^{s-1}dx  \leq\frac{C}{R^2}\int_{B_{2R}}v^{s}dx, \quad \forall\; R > 0.
\end{equation*}
\end{lem}
\noindent {\bf Proof}. Take $\phi \in C_0^2(\R^N).$  Let $(u,v)$ be a stable solution of
\eqref{1.1}, the integration by parts yields that
\begin{align}
\label{3.2}
\begin{split}
\int_{\mathbb{R}^N}|\nabla
u^{\frac{q+1}{2}}|^2\phi^2 dx & = \frac{(q+1)^2}{4}\int_{\mathbb{R}^N}u^{q-1}|\nabla u|^2\phi^2 dx\\ & = \frac{(q+1)^2}{4q}\int_{\mathbb{R}^N}\phi^2\nabla(u^q)\nabla
udx  \\
&= \frac{(q+1)^2}{4q}\int_{\mathbb{R}^N}\rho(x)u^qv^{P}\phi^2dx  + \frac{q+1}{4q}\int_{\mathbb{R}^N}u^{q+1}\Delta(\phi^2)dx ,
\end{split}
\end{align}
and
\begin{align}
\label{3.3a}
  (q+1)\int_{\mathbb{R}^N}u^q\phi\nabla u\nabla\phi dx  = \frac{1}{2}\int_{\mathbb{R}^N}\nabla(u^{q+1})\nabla(\phi^2)dx  = -\frac{1}{2}\int_{\mathbb{R}^N}u^{q+1}\D(\phi^2)dx .
\end{align}
 Take $\varphi =u^{\frac{q+1}{2}}\phi$ with $q > 0$ into the stability inequality \eqref{1.3} and using \eqref{3.2}-\eqref{3.3a}, we obtain
\begin{align*}
& \sqrt{p\theta}\int_{\mathbb{R}^N}\rho(x)u^{\frac{\theta-1}{2}}v^{\frac{p-1}{2}}u^{q+1}\phi^2dx\\ \leq &\; \int_{\mathbb{R}^N}|\nabla\varphi|^2dx\\
\leq & \;\frac{(q+1)^2}{4q}\int_{\mathbb{R}^N}\rho(x)u^qv^p\phi^2dx +C\int_{\mathbb{R}^N}u^{q+1}\Big[|\nabla\phi|^2+\Delta(\phi^2)\Big]dx,
\end{align*}
so we get
\begin{align*}
a_1\int_{\mathbb{R}^N}\rho(x)u^{\frac{\theta-1}{2}}v^{\frac{p-1}{2}}u^{q+1}\phi^2dx\leq
\int_{\mathbb{R}^N}\rho(x)u^qv^p\phi^2dx +C\int_{\mathbb{R}^N}u^{q+1}\Big[|\nabla\phi|^2+\Delta(\phi^2)\Big]dx,
\end{align*}
with $a_1=\frac{4q\sqrt{p\theta}}{(q+1)^2}$. Choose now
$\phi(x)= \varphi_0(R^{-1}x)$ where $\varphi_0 \in C_c^\infty(B_2)$ such that $\varphi_0
\equiv 1$ in $B_1$, there holds then
 \begin{equation}\label{3.3}
\int_{\mathbb{R}^N}\rho(x)u^{\frac{\theta-1}{2}}v^{\frac{p-1}{2}}u^{q+1}\phi^2dx\leq
\frac{1}{a_1}
\int_{\mathbb{R}^N}\rho(x)u^qv^p\phi^2dx+\frac{C}{R^2}\int_{B_{2R}}u^{q+1}dx.
\end{equation}
Similarly, applying the stability
inequality \eqref{1.3} with $\varphi = v^{\frac{r+1}{2}}\phi$, $r
> 0$, we obtain
\begin{equation}\label{3.4}
\int_{\mathbb{R}^N}\rho(x)u^{\frac{\theta-1}{2}}v^{\frac{p-1}{2}}v^{r+1}\phi^2dx\leq
\frac{1}{a_2} \int_{\mathbb{R}^N}\rho(x)u^\theta
v^r\phi^2dx+\frac{C}{R^2}\int_{B_{2R}}v^{r+1}dx
\end{equation}
 with $a_2=\frac{4r\sqrt{p\theta}}{(r+1)^2}$.
Combining \eqref{3.3} and \eqref{3.4},
 \begin{align}
 \label{3.5}
  \begin{split} & I_1+{a_2}^\frac{2(r+1)}{p+1} I_2\\
 := & \; \int_{\mathbb{R}^N}\rho(x)u^{\frac{\theta-1}{2}}v^{\frac{p-1}{2}}u^{q+1}\phi^2dx +
  {a_2}^\frac{2(r+1)}{p+1}\int_{\mathbb{R}^N}\rho(x)u^{\frac{\theta-1}{2}}v^{\frac{p-1}{2}}v^{r+1}\phi^2dx\\
  \leq & \; \frac{1}{a_1}
\int_{\mathbb{R}^N}\rho(x)u^qv^p\phi^2dx+{a_2}^\frac{2r+1-p}{p+1}\int_{\mathbb{R}^N}\rho(x)u^\theta
v^r\phi^2dx+\frac{C}{R^2}\int_{B_{2R}}\left(u^{q+1} +
v^{r+1}\right)dx.
  \end{split}
 \end{align}
Fix now
\begin{equation}\label{3.6}
 q=\frac{(\theta+1)r}{p+1}+\frac{\theta-p}{p+1}, \quad \mbox{ or equivalently }
q+1=\frac{(\theta+1)(r+1)}{p+1}.
\end{equation}
Let $r> \frac{p-1}{2},$ by Young's inequality, there holds
\begin{align*}
& \frac{1}{a_1}\int_{\mathbb{R}^N}\rho(x)u^qv^p\phi^2dx\\
= & \; \frac{1}{a_1}\int_{\mathbb{R}^N}\rho(x)u^{\frac{\theta-1}{2}}v^{\frac{p-1}{2}}u^{\frac{(\theta+1)r}{p+1}+\frac{\theta+1}{p+1}\left(\frac{1-p}{2}\right)}v^{\frac{p+1}{2}}\phi^2 dx \\
= & \;\frac{1}{a_1}\int_{\mathbb{R}^N}\rho(x)u^{\frac{\theta-1}{2}}v^{\frac{p-1}{2}}u^{(q+1)\frac{2r+1-p}{2(r+1)}}v^{\frac{p+1}{2}}\phi^2dx\\
\leq & \; \frac{2r+1-p}{2(r+1)}\int_{\mathbb{R}^N}\rho(x)u^{\frac{\theta-1}{2}}v^{\frac{p-1}{2}}u^{q+1}\phi^2dx
 +\frac{p+1}{2(r+1)}a_1^{-\frac{2(r+1)}{p+1}}\int_{\mathbb{R}^N}\rho(x)u^{\frac{\theta-1}{2}}v^{\frac{p-1}{2}}v^{r+1}\phi^2dx
\\
= & \; \frac{2r+1-p}{2(r+1)}I_1+\frac{p+1}{2(r+1)} a_1^{-\frac{2(r+1)}{p+1}} I_2;
\end{align*}
and similarly we have
 \begin{equation}\nonumber
  {a_2}^{\frac{2r+1-p}{p+1}}\int_{\mathbb{R}^N}\rho(x) u^\theta v^r\phi^2 dx \leq \frac{p+1}{2(r+1)} I_1
  + \frac{2r+1-p}{2(r+1)}{a_2}^\frac{2(r+1)}{p+1}I_2.
   \end{equation}
   Combining the above two estimates with \eqref{3.5}, we derive that
\begin{equation}
{a_2}^{\frac{2(r+1)}{p+1}}I_2\leq
\left[\frac{2r+1-p}{2(r+1)}{a_2}^{\frac{2(r+1)}{p+1}}+\frac{p+1}{2(r+1)}{a_1}^{\frac{-2(r+1)}{p+1}}\right]I_2 +
\frac{C}{R^2}\int_{B_{2R}}\left(u^{q+1}
+ v^{r+1}\right)dx, \nonumber
\end{equation}
hence
\begin{equation*}
\frac{p+1}{2(r+1)}\left[(a_1a_2)^{\frac{2(r+1)}{p+1}}-1\right]
I_2\leq CR^{-2}a_1^{\frac{2(r+1)}{p+1}}\int_{B_{2R}}\left(u^{q+1} +
v^{r+1}\right)dx.
\end{equation*}
Thus, if $a_1a_2 > 1$, by the choice of $\phi$,
\begin{equation*}
    \int_{B_R}\rho(x)u^{\frac{\theta-1}{2}}v^{\frac{p-1}{2}}v^{r+1} dx \leq I_2 \leq \frac{C}{R^2}\int_{B_{2R}}\left(u^{q+1} + v^{r+1}\right)dx .
\end{equation*}
Using \eqref{3.6} and \eqref{estS}, there hold $u^{q+1} \leq Cv^{r+1}$ and
$u^{\frac{\theta-1}{2}}v^{\frac{p-1}{2}}v^{r+1}\geq u^\theta v^r$.
Denote $s=r+1$, we conclude that if $a_1a_2
> 1$ and $s > \frac{p+1}{2}$,
\begin{equation*}
 \int_{B_R}\rho(x)u^\theta v^{s-1}dx  \leq
\frac{C}{R^2}\int_{B_{2R}}v^sdx .
\end{equation*}
Furthermore, we can check that $a_{1}a_2>1$ is equivalent to $L(s)<0$, the proof is completed. \qed

\medskip
We need also the following $L^1$ elliptic regularity result, see Lemma 5 in \cite{co}.
\begin{lem} \label{Al.2.2} For any $1\leq \beta<\frac{N}{N-2},$ there exists $C>0$ such that for any smooth non-negative function $w$, we have
\begin{equation*}
\left(\int_{B_{R_{k}}}w^{\beta} dx\right)^\frac{1}{\beta}\leq
CR^{N\big(\frac{1}{\beta}-1\big)+2}\int_{B_{R_{k+1}}}|w|dx+CR^{N\big(\frac{1}{\beta}-1\big)}\int_{B_{R_{k+1}}}wdx.
\end{equation*}
\end{lem}
Applying the above two lemmas, we establish the following result which plays an essential role in iteration process.
  \begin{lem}
\label{3.1}
Suppose that $\rho$ satisfies $(\star)$ and let $(u,v)$ be a classical stable solution of \eqref{1.1}, with $1< p \leq \theta$. Then for any  $1\leq \lambda<\frac{N}{N-2},$ $2t_0^-<q<s_0$ and
 nonnegative integer $k \geq 1$, there holds
\begin{equation}
\label{lnew2}
\left(\int_{B_{R_k}}v^{q\lambda}dx \right)^\frac{1}{\lambda}\leq
CR^{N\big(\frac{1}{\lambda}-1\big)}\int_{B_{R_{k+2}}}v^qdx, \quad \mbox{for all } R\geq 1.
\end{equation}
\end{lem}
\noindent {\bf Proof}. A simple calculation gives
\begin{equation*}
|\D (v^{q})|\leq q(q-1)v^{q-2}|\nabla v|^{2} +q\rho(x)v^{q-1}u^{\theta}.
\end{equation*}
 Using Lemma \ref{Al.2.2}, we get
\begin{align}\label{3.10a}
\begin{split}
\left(\int_{B_{R_{k}}}v^{q\lambda}dx \right)^\frac{1}{\lambda}\leq & \;
CR^{N\big(\frac{1}{\lambda}-1\big)+2}\int_{B_{R_{k+1}}}v^{q-2}|\nabla v|^{2}dx+CR^{N\big(\frac{1}{\lambda}-1\big)+2}\int_{B_{R_{k+1}}}\rho(x)v^{q-1}u^{\theta}dx\\ & \; + CR^{N\big(\frac{1}{\lambda}-1\big)}\int_{B_{R_{k+1}}}v^{q}dx.
\end{split}
\end{align}
Now, take a cut-off function $\phi\in C_{0}^{2}(B_{R_{k+2}})$ verifying $\phi\equiv 1$ in $B_{R_{k+1}}$ and $|\nabla\phi|\leq \frac{C}{R}.$ Multiplying $-\D v= \rho(x)u^{\theta}$ by $v^{q-1} \phi ^2$ and integrating by parts, we have
\begin{align}\label{3.9}
 (q-1)\int_{\mathbb{R}^N}v^{q-2}|\nabla v|^{2}\phi^{2}dx=-2\int_{\mathbb{R}^N}v^{q-1}\phi\nabla v\nabla \phi dx+\int_{\mathbb{R}^N}\rho(x)v^{q-1} u^{\theta}\phi ^2dx. \end{align}
By Young's inequality,
\begin{align*}
  2\int_{\mathbb{R}^N}v^{q-1}|\nabla v||\nabla \phi|\phi dx  \leq \frac{q-1}{2}\int_{\mathbb{R}^N}v^{q-2} |\nabla v|^{2}\phi ^2dx +
  C\int_{\mathbb{R}^N}v^{q}|\nabla \phi| ^2dx.
\end{align*}
Inserting this into \eqref{3.9}, using the properties of $\phi$, we obtain
\begin{align*}
\int_{\mathbb{R}^N}v^{q-2}|\nabla v|^{2}\phi^{2}dx \leq C\int_{B_{R_{k+2}}}\rho(x)v^{q-1}u^{\theta}dx+\frac{C}{R^{2}}\int_{B_{R_{k+2}}}v^{q}dx.
\end{align*}
Substituting the above inequality into \eqref{3.10a}, there holds
\begin{equation*}
\left(\int_{B_{R_{k}}}v^{q\lambda} dx\right)^\frac{1}{\lambda}\leq CR^{N\big(\frac{1}{\lambda}-1\big)+2}\int_{B_{R_{k+2}}}\rho(x)v^{q-1}u^{\theta}dx+CR^{N\big(\frac{1}{\lambda}-1\big)}\int_{B_{R_{k+2}}}v^{q}dx.
\end{equation*}
 Since $\rho$ satisfies $(\star),$ we can use Lemmas \ref{newl} to find \eqref{lnew2}.\qed

\medskip
 Now, we can follow exactly the iteration process as for Corollary 2 in \cite{co} (see also Proposition 3.1 in \cite{HU}) to obtain
\begin{cor}
\label{co3.1}
Suppose that $1< p \leq \theta$ and $\rho$ satisfies $(\star).$ Let $(u,v)$ be
a classical stable solution of (1.1) and $q \in (2t_0^-, s_0)$, then for $q \leq \beta<\frac{N}{N-2}s_0$,
there are $\ell\in \N$ and $C<\infty$ such that for any $R > 0$,
 \begin{align*}
 \left(\int_{B_R}v^{\beta} dx\right)^\frac{1}{\beta}\leq
CR^{N\big(\frac{1}{\beta}-\frac{1}{q}\big)}\left(\int_{B_{R_\ell}}v^q dx\right)^\frac{1}{q}, \quad \mbox{with } R_\ell = 2^\ell R.
\end{align*}
\end{cor}

Now we are in position to complete the proof of Theorem \ref{main1}.

\medskip\noindent
{\bf Proof of Theorem \ref{main1} completed.} Let $(u, v)$ be a classical stable solution of \eqref{1.1} with $\rho$ satisfying $(\star)$. We split the proof into two cases.

\medskip
\textit{Case 1: $p >\frac{4}{3}$.} Let $p>q>0.$ Using H\"older's inequality, there holds
   \begin{align}\label{3.13}
      \int_{B_{R}}v^qdx \leq\left(\int_{B_{R}}\rho v^pdx\right)^{\frac{q}{p}} \left(\int_{B_{R}}\rho ^{-\frac{p}{p-q}}dx\right)^{\frac{p-q}{p}}.
   \end{align}
   Applying Lemma \ref{l.2.1}, from $(\star)$ we get
    \begin{align*}
        \int_{B_{R}}v^qdx
         \leq C R^{\left[N-\frac{2(\theta+1)p}{p\theta-1}-\frac{(p+1)\alpha}{p\theta-1}\right]\frac{q}{p}
         +\left(N-\frac{\alpha q}{p-q}\right)\frac{p-q}{p}}
         = CR^{N-\frac{(2+\alpha)(\theta+1)}{p\theta-1}q}.
      \end{align*}
   By Remark \ref{rnew}, we known that $2t_0^-<p,$ then applying Corollary \ref{co3.1} with $2t_0^-<q<p$ and combining with Lemma \ref{l.2.1}, we can claim that for any $p \leq \beta < \frac{N}{N-2}s_0$, there
exists $C > 0$ such that
 \begin{align}
   \label{3.12}  \left(\int_{B_R}v^{\beta}dx \right)^\frac{1}{\beta}\leq
CR^{N\left(\frac{1}{\beta}-\frac{1}{q} \right)+\frac{1}{q}\left[N-\frac{(2+\alpha)(\theta+1)}{p\theta-1}q\right]}.
 \end{align}
Note that
$$N\left(\frac{1}{\beta}-\frac{1}{q} \right)+\frac{1}{q}\left[N-\frac{(2+\alpha)(\theta+1)}{p\theta-1}q\right] < 0 \quad \Leftrightarrow \quad N<\frac{(2+\alpha)(\theta+1)}{p\theta-1}\beta.$$
Suppose now
$$N<2+\left(\frac{(2+\alpha)(\theta+1)}{p\theta-1}\right)s_0.$$ We can take $\beta$ small but close to $\frac{N}{N-2}s_0$ such that $N<\frac{(2+\alpha)(\theta+1)}{p\theta-1}\beta$. With
a such $\beta$, tending $R\rightarrow\infty$ in \eqref{3.12}, we get $\|v\|_{L^\beta(\R^N)} = 0$, this is just impossible
since $v$ is positive. In other words, the equation \eqref{1.1} has no classical
stable solution if $N<2+(2+\alpha)x_0$ where $x_0=\frac{\theta+1}{p\theta-1}s_0$.

\medskip
Moreover, adopting the proof of Remark 2 in \cite{co1}, we can easily show that
\begin{align*}
2t_0^+\frac{\theta+1}{p\theta-1}>4,\quad \forall\,\theta\geq p>1.
  \end{align*}
By Lemma \ref{l.2.4}, $x_0 > 2t_0^+\frac{\theta+1}{p\theta-1} > 4$. This means that if $N\leq 10+4\alpha$, \eqref{1.1} has no classical stable solution
for any $\theta \geq p > \frac{4}{3}.$

\medskip
\textit{Case 2: $1<p\leq \frac{4}{3}$ and $u$ is bounded.} Let now $2>q>0,$ using \eqref{3.13}, with $p$ is replaced by $2$ and applying Lemma \ref{lemnew}, it follows that for any $R > 1$,
\begin{align*}
        \int_{B_{R}}v^qdx & \leq C R^{\left[N-\frac{2(\theta+1)p}{p\theta-1}-\frac{(p+1)\alpha}{p\theta-1} - \frac{2(2+\alpha)(2 - p)}{\theta+p-2}\right]\frac{q}{2}
         +\left(N-\frac{\alpha q}{2-q}\right)\frac{2-q}{2}}\\
         & = CR^{N-\left[\frac{(\theta+1)p}{p\theta-1} + \frac{(2+\alpha)(2 - p)}{\theta+p-2} +\frac{p(\theta+1)\alpha}{2(p\theta-1)}\right] q},
      \end{align*}
Proceeding as above, we can apply Corollary \ref{co3.1}, with $2t_0^-< q<2$ and $q<\beta<\frac{Ns_0}{N-2}$ to complete the proof of Theorem \ref{main1}. \qed

\bigskip\noindent
 \textbf{Proof of Corollary \ref{main2}.} Let $u$ be a solution of  the weighted Lane-Emden equation \eqref{D}, then $v = u$ verify the system \eqref{1.1} with $p = \theta$. Remark that $u$ is stable for \eqref{D} means just the estimate \eqref{1.3} holds true with $v=u$ and $p=\theta$, which is the departure point of our study. Moreover, we have $$t_0^{\pm}= p\pm\sqrt{p^2-p}$$ and
\begin{align*}
  L(s) = s^4 -16p^2s^2 +
32p^2s-16p^2=(s^2+4p(s-1))(s-2t_0^-)(s-2t_0^+).
\end{align*}
Then, $2t_0^+$ is the largest root of $L$ as $t_0^+ > p > 1$. Therefore
$$x_0 = \frac{2p+2\sqrt{p^2-p}}{p-1}$$ is the largest root of $H$, and we can check easily that $x_0 >4$ for all $p > 1.$ The result follows immediately by applying Theorem \ref{main1}. \qed

\bigskip\noindent
{\bf Acknowledgments}. We would like to thank Professor Dong Ye for suggesting us this problem and for many helpful comments.


\begin{thebibliography}{999}
\bibitem{cdg} W. Chen, L. Dupaigne and M. Ghergu, A new critical curve for the Lane-Emden system,
{\em Discrete Contin. Dyn. Syst.} \textbf{34} (2014), 2469-2479.

\bibitem{co} C. Cowan, Liouville theorems for stable Lane-Emden systems and biharmonic problems, {\em Nonlinearity} \textbf{26} (2013), 2357-2371.

\bibitem{co1}
C. Cowan, Regularity of stable solutions of a Lane-Emden type system, preprint (2012), arXiv:1204.4273.

\bibitem{cf}
C. Cowan and M. Fazly; On stable entire solutions of semilinear elliptic equations with weights, Proc. Amer.
Math. Soc. \textbf{140} (2012), 2003-2012

\bibitem{cg} C. Cowan and N. Ghoussoub, Regularity of semi-stable solutions to fourth order
nonlinear eigenvalue problems on general domains, {\em Calc. Var. PDE.} \textbf{49} (2014), 291-305.

\bibitem{ddww} J. D\'avila, L. Dupaigne, K. Wang and J. Wei, A monotonicity formula and a Liouville-type theorem
for a fourth order supercritical problem, {\em Adv. Math.} \textbf{258} (2014), 240-285.

\bibitem{dfs}
L. Dupaigne, A. Farina and B. Sirakov, Regularity of the extremal solutions for the Liouville system, in:
Geometric Partial Differential Equations, in: Publications of the Scuola Normale Superiore/CRM
Series, vol. \textbf{15} (2013), pp. 139-144.

\bibitem{Far} A. Farina, On the classification of solutions of the Lane-Emden equation on unbounded domains of $\R^N$, {\em J. Math. Pures Appl.} \textbf{87} (2007), 537-561.

\bibitem{Fa}
M. Fazly, Liouville type theorems for stable solutions of certain elliptic systems, \emph{Adv. Nonlinear Stud.} \textbf{12} (2012), 1-17.

\bibitem{gnw} C. Gui, W. Ni and X. Wang, On the stability and instability of positive steady states of a semilinear heat equation in ${\bf R}^n$, {\em Comm. Pure Appl. Math.} \textbf{Vol. XLV} (1992), 1153-1181.

\bibitem{hhy}
H. Hajlaoui, A. Harrabi and D. Ye, On stable solutions of the biharmonic problem with polynomial growth, {\em Pacific J. Math.} \textbf{270} (2014), 79-93.

\bibitem{HU}
L. Hu, Liouville type results for semi-stable solutions of the weighted Lane-Emden
system, \emph{J. Math. Anal. Appl.} \textbf{432} (2015), 429-440.

\bibitem{jl}
D.D. Joseph and T.S. Lundgren, Quasilinear Dirichlet problems driven
by positive sources, \emph{Arch. Rational Mech. Anal.} {\bf 49}
(1973), 241-269.

\bibitem{MPO}
E. Mitidieri and S. Pohozaev, A priori estimates and the absence of solutions of nonlinear
partial differential equations and inequalities, \emph{Tr. Mat. Inst. Steklova} \textbf{234} (2001), 1-384.

\bibitem{Mo}
M. Montenegro, Minimal solutions for a class of elliptic systems, \emph{Bull. London Math. Soc.} \textbf{37} (2005),
405-416.

\bibitem{sz} J. Serrin and H. Zou, Non-existence of positive solutions
of Lane-Emden systems, {\em Diff. Inte. Equations} \textbf{9} (1996),
635-653.

\bibitem{S} P. Souplet, The proof of the Lane-Emden conjecture in four space dimensions, {\em Adv.  Math.} \textbf{221} (2009), 1409-1427.

\bibitem{wy} J. Wei and D. Ye, Liouville theorems for stable solutions of biharmonic problem, {\em Math. Ann.} \textbf{356} (2013), 1599-1612.

\end{thebibliography}
\end{document}